\documentclass[a4paper,11pt]{amsart}
\usepackage[utf8]{inputenc}
\usepackage[T1]{fontenc}

\usepackage{amsmath}

\usepackage[all]{xy}
\usepackage{amssymb}
\usepackage{fullpage}
\usepackage{color}
\usepackage{stmaryrd}
\usepackage{tikz-cd}

\newcommand{\N}{\mathbb{N}}
\newcommand{\Q}{\mathbb{Q}}
\newcommand{\Z}{\mathbb{Z}}
\newcommand{\R}{\mathbb{R}}
\newcommand{\C}{\mathbb{C}}

\newcommand{\rH}{\mathrm{H}}
\newcommand{\Hom}{\mathrm{Hom}}
\newcommand{\Hit}{\mathrm{Hit}}

\newcommand{\Aut}{\mathrm{Aut}}
\newcommand{\Int}{\mathrm{Int}}
\newcommand{\Out}{\mathrm{Out}}

\newcommand{\cE}{\mathcal{E}}
\newcommand{\cM}{\mathcal{M}}

\newcommand{\cC}{\mathcal{C}}

\newcommand{\cK}{\mathcal{K}}
\newcommand{\cT}{\mathcal{T}}

\newcommand{\frg}{\mathfrak{g}}
\newcommand{\fG}{\mathcal{G}}

\newcommand{\fX}{\mathcal{X}}
\newcommand{\fL}{\mathfrak{L}}

\newcommand{\bH}{\mathbf{H}}

\newcommand{\PGL}{\mathbf{PGL}}
\newcommand{\GL}{\mathbf{GL}}
\newcommand{\PO}{\mathbf{PO}}
\newcommand{\PSp}{\mathbf{PSp}}

\newcommand{\si}{\sigma}

\newcommand{\lra}{\longrightarrow}
\newcommand{\eps}{\varepsilon}
\newcommand{\ga}{\gamma}
\newcommand{\Ga}{\mathrm{\Gamma}}

\newcommand{\lmt}{\longmapsto}

\newcommand{\piY}{\pi_1Y}
\newcommand{\piX}{\pi_1X}

\newcommand{\Fix}{\mathrm{Fix}}

\newcommand{\pifX}{\pi_1\fX}

\renewcommand{\rho}{\varrho}
\renewcommand{\phi}{\varphi}
\renewcommand{\geq}{\geqslant}
\renewcommand{\leq}{\leqslant}

\numberwithin{equation}{section}

\theoremstyle{definition}
\newtheorem{definition}{Definition}[section]

\theoremstyle{plain}
\newtheorem{theorem}[definition]{Theorem}

\allowdisplaybreaks[4]

\title{Symmetric differentials and the dimension of Hitchin components for orbi-curves}
\author{Florent Schaffhauser}
\address{Universidad de Los Andes, Departamento de Matem\'aticas, Carrera 1 \#18A-12, 111 711 Bogot\'a, Colombia\\ \& Universit\'e de Strasbourg, IRMA UMR 7501, 7 rue Ren\'e Descartes, 67 000 Strasbourg, France.}
\email{\texttt{schaffhauser@math.unistra.fr}}

\begin{document}

\begin{abstract}This note is based on a talk given at the 2019 ISAAC Congress in Aveiro. We give an expository account of joint work with Daniele Alessandrini and Gye-Seon Lee on Hitchin components for orbifold groups, recasting part of it in the language of analytic orbi-curves. This reduces the computation of the dimension of the Hitchin component for orbifold groups to an application of the orbifold Riemann-Roch theorem.
\end{abstract}

\maketitle

\section{Hitchin components for orbifold fundamental groups}

\subsection{Compact orbi-surfaces of negative Euler characteristic}

An \emph{orbifold} is a kind of space that generalises the notion of a manifold (be it a topological, differentiable or analytic one). For instance, a differentiable orbifold is a type of space that locally looks like the quotient of an open set $U\subset \R^n$ by a finite group of diffeomorphisms $\Ga\subset \mathrm{Diff}(U)$. What is meant here by \emph{quotient} depends a lot on how one understands the expression \emph{a type of space}. For us, it will be sufficient to consider (topological, differentiable or analytic) \emph{stacks} as our notion of space. Such a stack is then called an orbifold if it admits a \emph{covering} by open substacks of the form $[U/\Ga]$, parameterising families of $\Ga$-orbits in $U$, where $U$ is the local model for representable stacks (i.e.\ manifolds) and $\Ga$ is a finite subgroup of the automorphism group of $U$. A fundamental example of orbifold is the stack $\mathcal{X}:=[M/\pi]$, where $\pi$ is a discrete group acting (effectively and) properly on a manifold $M$.

A \textit{coarse moduli space} (CMS) for an orbifold $\mathcal{X}$ is a manifold $X$ equipped with a morphism $p:\fX\lra X$ that satisfies the following universal property for all manifolds $M$:
$$
\xymatrix{
\mathcal{X} \ar[r] \ar[d]^p & M\\
X \ar@{-->}[ur]^{\exists!}
}
$$ In the first part of the paper, we will work only with (effective) \emph{differentiable} orbifolds. Then, up to real dimension $2$, it suffices to enlarge the category of manifolds slightly and accommodate manifolds with corners, to ensure that coarse moduli spaces always exist. This is convenient because it allows us to think of an orbi-surface (or even an orbi-surface with boundary) as an ordinary surface with an extra structure, namely some ``special points", all of whose open neighbourhoods are of the form $[U/\Ga]$ with non-trivial $\Ga$. As a matter of fact, since we are in the $\cC^\infty$ setting and $\Ga$ is finite, we can always assume that it acts on the open set $U\subset \R^2$ preserving a positive-definite metric. The classification of linear isometries of the Euclidean plane then tells that a point in the coarse moduli space $U/\Ga$ is of one of the following three types:
\begin{enumerate}
\item A \emph{cone point}, which admits an open neighbourhood of the form $D(0;\eps) / \mathrm{C}_m$, where $\mathrm{C}_m\simeq\Z/m\Z$ is a finite cyclic group of order $m$, acting on the open disk $D(0;\eps)$ by rotation. Such a cone point is said to have \emph{order} $m$.
\item A \emph{dihedral point} (also called \emph{corner reflector}), which admits an open neighbourhood of the form $D(0;\eps) / \mathrm{D}_m$, where $\mathrm{D}_m\simeq \mathrm{C}_m\rtimes \Z/2\Z$ is the dihedral group of order $2m$. Such a dihedral point is said to have \emph{order} $m$.
\item A \emph{mirror point}, which admits an open neighbourhood of the from $D(0;\eps) / \Z/2\Z$, where $\Z/2\Z$ acts on $D(0;\eps)$ by reflection through a diameter.
\end{enumerate}

For instance, a $2$-dimensional orbifold could have a triangle for a coarse moduli space: the edges are mirror points, while the vertices are dihedral points. Another good thing about (compact) orbi-surfaces is that they admit an \textit{orbifold Euler characteristic}, computable explicitly from the coarse moduli space through the following formula (in which $k$ is the number of cone points, $\ell$ the number of dihedral points, $m_i$ is the order of the $i$-th cone point and $n_j$ is the order of the $j$-th dihedral point): \begin{equation}\label{Euler_char} \chi(\mathcal{X}) = \chi(X) - \textstyle\sum_{i=1}^{k} \left(1-\frac{1}{m_i}\right) - \textstyle\frac{1}{2}  \textstyle\sum_{j=1}^{\ell} \left(1-\frac{1}{n_j}\right)\in\Q.\end{equation} In the so-called orientable case, the quantity $\chi(\mathcal{X})$ is negative if and only if the coarse moduli space $X$ is a closed surface of genus at least $2$, or a torus with at least one cone point, or a sphere with at least three cone points. For a more complete treatment of the fundamental properties of orbifolds, we refer for instance to \cite{Thurston,Scott,Cooper_et_al,CG} and for the stacky point of view, we refer to \cite{Hinich_Vaintrob,BGNX_string_top}. An important property of the Euler characteristic is that it is multiplicative: if $\fX\simeq [\mathcal{Y}/\pi]$, then $\chi(\fX)=\frac{\chi(\mathcal{Y})}{|\pi|}$.

\subsection{Fundamental group and hyperbolic structures}

A \emph{cover} of an orbifold $\mathcal{X}$ is a morphism $\mathcal{Y} \lra \mathcal{X}$ which, in orbifold charts, is conjugate to a morphism of the form $\amalg_{i\in I}[U/\Ga_i]\lra[U/\Ga]$, where each $\Ga_i$ is a subgroup of $\Ga$. In particular, the canonical map $U\lra [U/\Ga]$ is an orbifold cover. A more concrete example is given as follows: the ``flattening" of a sphere with $3$ cone points is a $2$-to-$1$ cover of a triangle, the cone points upstairs being mapped to dihedral points of the same order downstairs. There is an orbifold structure on the inverse limit of all connected covers and the latter is called the \emph{universal cover} of $\mathcal{X}$. The fundamental group of $\mathcal{X}$ is the automorphism group of the universal cover (whose total space may or may not be a manifold). We denote by $\pi_1\mathcal{X}$ the fundamental group of $\mathcal{X}$. For instance, if $M$ is a simply connected manifold and $\pi$ is a discrete group acting properly on $M$, then $\pi_1([M/\pi])\simeq \pi$. If one is careful about base points, connected covers of an orbifold $\mathcal{X}$ correspond bijectively to subgroups of $\pi_1\mathcal{X}$.

A \emph{hyperbolic structure} on a differentiable orbifold $\mathcal{X}$ is a covering by open substacks of the form $[U/\Ga]$ in which $U\subset \bH^2$ is an open subspace \emph{of the (real) hyperbolic plane} and $\Ga\subset\mathrm{Isom}(\bH^2)\simeq\PGL(2;\R)$ is a finite subgroup of the isometry group of $\bH^2$ that leaves $U$ invariant. If $\chi(\mathcal{X})<0$, then $\mathcal{X}$ admits hyperbolic structures and its universal cover is isomorphic to $\bH^2$. The deformation space of hyperbolic structures on $\mathcal{X}$ is identified, via the space of holonomy representations of such structures, to a connected component of the topological space $$\Hom\big(\pi_1\mathcal{X};\PGL(2;\R)\big)/\, \PGL(2;\R).$$ Namely, it is the space of \emph{discrete and faithful} representations $\rho:\pi_1\mathcal{X}\lra\PGL(2;\R)$. Thus, if $\chi(\mathcal{X})<0$, it is always possible to identify $\pi_1\mathcal{X}$ with a discrete subgroup of $\PGL(2;\R)$, i.e.\ a \emph{Fuchsian group}. If the orbifold $\mathcal{X}$ is orientable (which in dimension $2$ amounts to saying that the CMS $X$ is an orientable surface and that all groups $\Ga$ appearing in the orbifold charts contain only orientation-preserving transformations), then the fundamental group of $\fX$ admits the following presentation: 
\begin{equation}\label{fundamental_group_cx_orbicurve}
\pi_1\mathcal{X} \simeq \left<(a_i,b_i)_{1\leq i\leq g}, (c_j)_{1\leq j\leq k}\ |\ \textstyle{\prod_{1\leq i \leq g}[a_i,b_i] \prod_{1\leq j\leq k} c_j}=1= c_1^{m_1}=\ldots=c_k^{m_k} \right>.
\end{equation} We will denote by $\pi_{g,(m_1,\,\ldots\,,m_k)}$ the group defined by the presentation \eqref{fundamental_group_cx_orbicurve}. The space of discrete and faithful representations of $\pi_1\mathcal{X}$ in $\PGL(2;\R)$ will be called the \emph{Teichmüller space} of $\mathcal{X}$ and denoted by $\cT(\fX)$. It is homeomorphic to a real vector space of dimension $-3\chi(X) + 2k+\ell$. In particular, it is reduced to a point if $\mathcal{X}$ is a (quotient of a) sphere with three cone points (see for instance \cite{CG} for a full account on this, including for the more refined notion of orbifold with boundary).

\subsection{Hitchin components}

Let $\frg_\C$ be a simple complex Lie algebra. The \emph{adjoint group} $G_\C:=\Int(\frg_\C)$ is the neutral component, in the Lie group topology, of $\Aut(\frg_\C)$, and it is a complex Lie group with trivial centre, whose Lie algebra is isomorphic to $\frg_\C$. Given a real form $\frg$ of $\frg_\C$, there is an associated anti-holomorphic involution $\theta$ of $G_\C$, whose fixed-point set we denote by $G$. It consists of interior automorphisms of $\frg$ that commute with $\theta$. The neutral component of $G$ is $\Int(\frg)$. In particular, $G$ is not necessarily connected. For instance, if we choose $\frg=\mathfrak{sl}(n;\R)$, then $G\simeq\PGL(n;\R)$, which is connected if $n$ is odd and has two connected components if $n$ is even. In what follows, we shall always assume that $\frg$ is the \emph{split real form} of the Lie algebra $\frg_\C$.

In \cite{Hitchin_Teich}, N.~Hitchin studies representations of surface groups into $G$ and shows that the representation space $\Hom(\piX;G)/G$ has a contractible connected component. His definition of that component rests on the notion of \emph{Fuchsian representation}, which itself depends on the choice of a so-called \emph{principal morphism} $\kappa:\PGL(2;\R)\lra G$, first introduced by B.~Kostant (\cite{Kostant_ppal_1959}). When $G=\PGL(n;\R)$, this morphism is induced by the linear action of $\GL(2;\R)$ on the space $V_n$ of homogeneous polynomials of degree $n-1$ in two variables $x,y$. Hitchin's definition, extended to the orbifold case, is then the following.
Given an orbi-surface $\fX$ of negative Euler characteristic and a principal morphism $\kappa$ from $\PGL(2;\R)$ to the split real form $G$ of $\Int(\frg_\C)$, a representation $\rho:\pifX\lra G$ is called \emph{Fuchsian} if it lifts to a discrete and faithful representation $h:\pifX\lra \PGL(2;\R)$, in the sense that the following diagram becomes commutative:
$$
\xymatrix{
 & \PGL(2;\R) 
\ar^{\kappa}[d] \\
\pifX \ar@{-->}[ur]^{h}
\ar[r]^{\widehat{\rho}} & G
}
$$ This defines a map $\cT(\fX)\lra \Hom(\pifX;G)/G$ whose image is called the \emph{Fuchsian locus}. As the Teichmüller space $\cT(\fX)$ is connected, this map picks out a single connected component of the representation space $\Hom(\pifX;G)/G$, called the \emph{Hitchin component} and denoted by $\Hit(\pifX;G)$. When $G=\PGL(2;\R)$, we have $\Hit(\pifX;\PGL(2;R))\simeq\cT(\fX)$, by definition. For split real groups $G$ of \emph{higher rank}, Hitchin components form a family of so-called \emph{Higher Teichmüller spaces} (\cite{Wienhard_ICM}). Indeed, Hitchin representations are discrete and faithful (\cite{Labourie_Anosov,Fock_Goncharov_IHES}). In the surface group case, Hitchin has proved that $\Hit(\piX,G)$ has a trivial topology:

\begin{theorem}[Hitchin, \cite{Hitchin_Teich}]\label{dim_HC_surface_case}
Let $X$ be a closed orientable surface of negative Euler characteristic and let $G$ be the split real form of $\Int(\frg_\C)$, where $\frg_\C$ is a simple complex Lie algebra. Then $\Hit(\piX;G)$ is homeomorphic to a real vector space of dimension $-\chi(X)\dim G$.
\end{theorem}

This formula cannot be generalised directly to the orbifold case, as $\chi(\fX)$ is not an integer in general. However, for $G=\PGL(3;\R)$, S.~Choi and W.~Goldman have proved the following formula.

\begin{theorem}[Choi \& Goldman, \cite{CG}]\label{CG_thm}
Let $\fX$ be a closed orbi-surface of negative Euler characteristic and with coarse moduli space $X$. Then $\Hit(\pifX;\PGL(3;\R))$  is homeomorphic to a real vector space of dimension $-8\chi(X) + (6k-2k_2) + (3\ell-\ell_2)$, where $k_2$ (respectively, $\ell_2$) is the number of cone points (respectively, dihedral points) of order $2$ of $\fX$.
\end{theorem}

In collaboration with D.~Alessandrini and G.S.~Lee, we have been looking at Hitchin components for orbifold groups and we have obtained the following common generalisation of the two results above. For the sake of clarity, we will present it here for the group $G=\PGL(n;\R)$ only, but our results hold for split real forms of all adjoint groups of simple complex Lie algebras, for instance $\PSp^\pm(2m;\R)$, $\PO(m,m+1)$, $\PO^\pm(m,m)$, or the exceptional Lie group $\mathbf{G}_2$.

\begin{theorem}[\cite{ALS}]\label{dim_HC_orb_case}
Let $\fX$ be a closed orbi-surface of negative Euler characteristic and with coarse moduli space $X$. Then $\Hit(\pifX;\PGL(n;\R))$ is homeomorphic to a real vector space of dimension $$-(n^2-1)\chi(X) + \textstyle\sum_{d=2}^n \left( 2 \textstyle\sum_{i=1}^k R(d,m_i) + \textstyle\sum_{j=1}^{\ell} R(d,n_j)\right)$$ where $m_i$ (respectively, $n_j)$ is the order of the $i$-th cone point (respectively, the $j$-th dihedral point) of $\fX$ and $R(d,m):=\left\lfloor d - \frac{d}{m}\right\rfloor$ is the integral part of the real number $\left(d-\frac{d}{m}\right)$.
\end{theorem}

As a matter of fact, like Choi and Goldman in \cite{CG}, we can also deal with the case of orbifolds with boundary. We also note that, when $Y$ has at most mirror points as orbifold singularities (no cone or dihedral points), then $\chi(\fX)=\chi(X)$ and Hitchin's formula holds without modifications. There is another way of writing the formula in Theorem \ref{dim_HC_orb_case}, which resembles more that of Theorem \ref{CG_thm}, and we refer to \cite{ALS} for it. We get for instance $$\dim\Hit\big(\pifX;\PGL(4;\R)\big) = -15\chi(X) + (12k - 4k_2 -2 k_3) + (6\ell - 2\ell_2 -\ell_3),$$ where again $k_i$ (respectively, $\ell_i$) is the number of cone points (respectively, dihedral points) of order $i$ of $\fX$. We see that this dimension may vanish for certain orbifolds $\fX$ and that such orbifolds form an infinite family, containing for instance all spheres with three cone points of order $(2,3,r)$ for all $r\geq 7$. This has applications to the rigidity of projective structures on Seifert-fibered spaces with base $\fX$ (see \cite{ALS} for details).

The methods of proof for Theorems \ref{dim_HC_surface_case} and \ref{CG_thm} are quite different. Hitchin uses tools from analytic and differential geometry (namely, Higgs bundles and the Non-Abelian Hodge Correspondence), while Choi and Goldman's methods are based on the interpretation of $\Hit(\pi\fX;G)$ as the deformation space of \emph{convex projective structures} on $\fX$. In the absence of such a geometric interpretation for general $G$, our approach in \cite{ALS} consists in adapting Hitchin's method to our setting. Thanks to an orbifold version of the Non-Abelian Hodge Correspondence, we show that the Hitchin component is homeomorphic to a space of symmetric differentials on an analytic orbi-curve, the dimension of which we can compute using the orbifold Riemann-Roch theorem, similarly to Hitchin's proof in the surface group case. We explain this in greater detail in the next section.

\section{Analytic parameterisation of Hitchin components}

\subsection{Analytic orbi-curves}

A \emph{complex analytic orbifold} is an analytic stack $\fX$ (over complex analytic manifolds) that admits a covering by open substacks of the form $[U/\Ga]$, where $U\subset \C^n$ is an open subset and $\Ga\subset \Aut(U)$ is a finite group of holomorphic transformations of $U$. If the open sets $U$ are all of complex dimension $1$, we say that $\fX$ is an \emph{orbi-curve} or an \emph{orbi-Riemann surface}. The only possible orbifold points in this case are cone points and it is a remarkable fact that there always exist coarse moduli spaces: an orbi-Riemann surface always has an ``underlying" Riemann surface, because if $\Ga\simeq \mathrm{C}_m$ acts by rotation of angle $\frac{2\pi}{m}$ on the open disk $D(0;\eps)$ then the map $z\lmt z^m$ induces a holomorphic chart $D(0;\eps)/\mathrm{C}_m \simeq D(0;\eps^m)$. In fact, the whole theory of complex analytic orbi-curves can be phrased in terms of \emph{Riemann surfaces with signature}, where the signature is the map $X\lra \N$ taking a point to its order (so the map is constant equal to $1$, except possibly over a finite set of points in $X$). We prefer to work, however, in the orbifold setting. In particular, subgroups of the orbifold fundamental group \eqref{fundamental_group_cx_orbicurve} correspond to connected analytic covers of the compact orbi-curve $\fX:=[\bH^2/\pi_{g,(m_1,\,\ldots\, , m_k)}]$.

To prove Theorem \ref{dim_HC_orb_case}, complex analytic orbi-curves will not be quite enough if we want to include the case of non-orientable differentiable orbi-surfaces. To deal with those, we need to consider also orbi-curves which are \emph{defined over the real numbers}. This essentially means complex analytic orbi-curves $\fX^+$ equipped with an anti-analytic involution $\si:\fX^+\lra\fX^+$ given, in local charts, by a $\Ga$-equivariant anti-holomorphic involution $\si:U\lra U'$. In particular, the orders of the points $x$ and $\si(x)$ have to coincide for all $x$. More intrinsically perhaps, one could consider \emph{dianalytic orbifolds}, for which local models are quotient stacks $[U/\Ga]$, where $U\subset \C^n$ is an open subset but the finite group $\Ga\subset \Aut^\pm(U)$ is now allowed to also contain anti-holomorphic transformations of $U$. If we consider such a dianalytic orbifold $\fX$, its fundamental group $\pi:=\pifX$ has a subgroup $\pi^+$, of index at most $2$, consisting of transformations that preserve the orientation of the universal cover $\widetilde{\fX}$, the latter being necessarily complex analytic: the quotient orbifold $\fX^+:= [\widetilde{\fX}/\pi^+]$ is a complex analytic orbifold which is a cover, of degree at most $2$, of $\fX$. If $\pi^+\neq\pi$, then $\pi/\pi^+$ acts on $\fX^+$ via an anti-holomorphic involution $\si$ and $\fX\simeq [\fX^+/\left<\si\right>]$. In this case, there is a short exact sequence $$1\lra \pifX^+ \lra \pifX \lra \{\pm 1\} \lra 1.$$ Note that $\fX^+$ has two cone points $x$ and $\si(x)$ (of the same order) for each cone point of $\fX$, and one cone point which is fixed by $\si$ for each dihedral point of $\fX$. Consider for instance the fundamental group of a triangle $\fX$ with vertices of respective orders $(p,q,r)$. The double cover $\fX^+$ is a sphere with three cone points, of respective orders $(p,q,r)$. The fundamental group of $\fX^+$ is the Von Dyck group $$\pi_{0,(p,q,r)}\simeq\left<a,b,c\ |\ a^p=b^q=c^r=abc=1\right>$$ of \eqref{fundamental_group_cx_orbicurve}, while that of $\fX$ is the Coxeter (triangle) group with presentation 
\begin{equation}\label{Coxeter_gp}
T_{(p,q,r)}:=\left<x,y,z\ |\ x^2=y^2=z^2=(xy)^p=(yz)^q=(zx)^r=1\right>.
\end{equation} The covering (flattening map) $\fX^+\lra\fX$ induces the injective group morphism $$\pi_{0,(p,q,r)}\lra T_{(p,q,r)}$$ defined by $a\lmt xy$, $b\lmt yz$, $c\lmt zx$, and the quotient map $T_{(p,q,r)} \lra \{\pm1\}$ is given by the reduced word length modulo $2$.

When $\fX$ is a compact orbi-curve of negative Euler characteristic, the fundamental group $\pifX$ is a finitely generated group that embeds onto a discrete subgroup of $\PGL(2;\R)$. Therefore, by Selberg's lemma, it contains a finite index normal subgroup which is torsion-free (\cite{Selberg}). Geometrically, this means that there exists a compact Riemann surface $Y$ and a finite Galois cover $Y\lra\fX$. If we denote by $\pi$ the automorphism group of that cover, we therefore have an isomorphism of orbifolds $[Y/\pi]\simeq\fX$, and a short exact sequence $$1\lra\piY\lra\pifX\lra\pi\lra 1.$$

\subsection{The Riemann-Roch formula}\label{RR}

An \emph{orbifold line bundle} $\fL$ over $\fX$ is a morphism of stacks $\fL \lra  \fX$ which is locally conjugate, in the orbifold chart $[U/\Ga]$ about $x$, to the orbifold $[(U\times \C)/\Ga]$, where $\Ga$ acts on $U\times \C$ via a linear representation $\rho_\Ga:\Ga\lra\GL(1,\C)$. When the finite group $\Ga$ is cyclic of order $m$, the morphism $\rho_\Ga$ sends a generator of $\Ga$ to an $m$-th root of unity. If we choose a generator $\ga$ of $\Ga$ and a primitive $m$-th root of unity $\zeta$, then $\rho_\Ga(\ga)=\zeta^a$ for a certain $a\in\{0;\,\ldots\, ;m-1\}$ which does not depend on the choices just made and is sometimes called the \emph{isotropy} at the point $x$. When $\Ga$ is a dihedral group, we write $\Ga\simeq \mathrm{C}_m\rtimes\Z/2\Z$, where $\Z/2\Z$ acts on the cyclic group $\mathrm{C}_m$ by inversion, and think of $\rho_\Ga:\Ga\lra\C^*$ as a morphism $\rho_{\mathrm{C}_m}:\mathrm{C}_m\lra \C^*$ as before, which in addition is $\Z/2\Z$-equivariant with respect to complex conjugation on $\C^*$. In particular, the number $a\in\{0;\,\ldots\, ;m-1\}$ again completely determines the morphism $\rho_\Ga$. Given a cone or dihedral point $x$ of order $m$, the quantity $\frac{a}{m}$, where $a\in\{0;\,\ldots\, ;m-1\}$ is defined as above, will be called the \emph{age} of the orbifold line bundle $\fL$ at $x$ and denoted by $\mathrm{age}_{x}(\fL)$. Consider for instance the canonical line bundle $K_\fX$ of an analytic orbi-curve $\fX$. The age of the tangent bundle at a cone point of order $m$ is $\frac{1}{m}$ (the action of $\mathrm{C}_m$ on tangent vectors being multiplication by a primitive root of unity) and, since the  canonical bundle is the dual of the tangent bundle in this case, the group $\Ga\simeq \mathrm{C}_m$ acts on tangent covectors at a point via multiplication by $\zeta^{-1}$, so the age of $K_\fX$ at a cone point of order $m$ is $\frac{m-1}{m}$. If we now look at tensor powers of $K_\fX$, then the action of $\mathrm{C}_m$ on homogeneous polynomial functions of degree $d$ over the tangent space at a cone point is given by multiplication by $\zeta^{d(m-1)}$, so the age of $K_\fX^d$ at a cone point is $$\textstyle\frac{d(m-1)\,\mathrm{mod}\, m}{m} = \textstyle\frac{d(m-1)}{m}-\left\lfloor \textstyle\frac{d(m-1)}{m}\right\rfloor.$$ We will see in Section \ref{computation_of_the_dimension} below that this is the origin of the term $R(d,m)=\left\lfloor \frac{d(m-1)}{m}\right\rfloor$ in Theorem \ref{dim_HC_orb_case}.

Let us denote by $\underline{\fL}$ the sheaf of local sections of $\fL$. There are associated cohomology groups $\rH^0(\fX;\underline{\fL})$ and $\rH^1(\fX;\underline{\fL})$, which are finite-dimensional complex or real vector spaces (depending on the field of definition of $\fX$). The \emph{Euler characteristic} of $\underline{\fL}$ is the integer $\chi(\fX;\underline{\fL}):= \dim\rH^0(\fX;\underline{\fL}) - \dim\rH^1(\fX;\underline{\fL})$. The Riemann-Roch formula computes this quantity by comparing it to the Euler characteristic of the structure sheaf $O_\fX$. To state the result, we still need the notion of \emph{degree} of an orbifold line bundle, of which we recall the following two definitions (in the complex case). When $\fX\simeq[Y/\pi]$, where $Y$ is a compact Riemann surface and $\pi$ is a finite group of analytic transformations of $Y$, an orbifold line bundle $\fL\lra\fX$ pulls back to a $\Ga$-equivariant analytic line bundle $\cE\lra Y$ and we can define the degree of $\fL$ as $\frac{\deg(\cE)}{|\pi|}\in\Q$, since this quantity is independent of the choice of the finite Galois cover $Y\lra\fX$. Equivalently, if we denote by $p:\fX\lra X$ the coarse moduli space of $\fX$, then, given an orbifold line bundle $\fL\lra\fX$, there exists a unique analytic line bundle $L\lra X$ and for each cone point $x_i$ of $\fX$ a well-defined integer $a_i\in\{0;\,\ldots\, ;m_i-1\}$ such that $$\fL\ \simeq\ p^*L\otimes O_{\fX}\big(\textstyle\sum_{i=1}^k a_ix_i\big).$$ We then have $\mathrm{age}_{x_i}(\fL) = \frac{a_i}{m_i}$ and $\deg(\fL):= \deg(L)+\sum_{i=1}^k \frac{a_i}{m_i}$, where $m_i$ is the order of the cone point $x_i$. For instance, when $\fL=K_{\fX}^d\simeq [K_Y^d/\pi]$, one can check that $$K_{\fX}^d\ \simeq\ p^*\big[K_X^d\otimes O_X\big(\textstyle\sum_{i=1}^k R(d,m_i)p(x_i)\big)\big] \otimes O_{\fX}\big(\textstyle\sum_{i=1}^k (d(m_i-1)\,\mathrm{mod}\,m_i) x_i\big)$$ so, using \eqref{Euler_char}, 
\begin{equation}\label{degree_powers_of_can_bdle}
\deg(K_{\fX}^d)= d(2g-2) + \textstyle\sum_{i=1}^k R(d,m_i) + \textstyle\sum_{i=1}^k \textstyle\frac{d(m_i-1)\,\mathrm{mod}\,m_i}{m_i} =-d\chi(\fX)=-d\textstyle\frac{\chi(Y)}{|\pi|} = \textstyle\frac{\deg(K_Y^d)}{|\pi|}\, ,
\end{equation} where $g:=\dim\rH^1(\fX;O_{\fX})$ is the genus of $\fX$. Indeed, $\pi$-invariant holomorphic sections of $K_Y^d$ correspond bijectively to meromorphic sections of $K_X^d$ with poles of order at most $R(d,m_i)$ at $x_i$, for all $i\in\{1;\,\ldots\,;k\}$. More generally, for $j=0,1$, there are isomorphisms $\rH^j(\fX;\underline{\fL}) \simeq \Fix_{\pi}\, \rH^j(Y;\underline{\cE}) \simeq \rH^j(X;\underline{L})$, from which one can deduce the following orbifold Riemann-Roch formula (see for instance \cite{Nasatyr_Steer_HYM}, \cite{Abramovic} or \cite{Liu_orbi_RR} for an exposition; the theorem itself is due to Kawasaki, \cite{Kawasaki}).

\begin{theorem}[Orbifold Riemann-Roch, \cite{Kawasaki}]\label{RR_thm}
Let $\fX$ be a compact complex analytic orbi-curve and denote its cone points by $(x_i)_{1\leq i\leq k}$. Let $\fL$ be an analytic line bundle over $\fX$. Then $$\chi(\fX;\underline{\fL})\ =\ \chi(\fX;O_{\fX}) + \deg\fL - \textstyle\sum_{i=1}^k \mathrm{age}_{x_i}(\fL).$$
\end{theorem}

For instance, $\chi(\fX;\underline{K_{\fX}}) = g-1$ and $\chi(\fX;\underline{K_{\fX}^2}) =3(g-1) + k$. When $\fX$ is defined over $\R$, we can deduce the appropriate version of the Riemann-Roch formula from the complex case, by applying Theorem \ref{RR_thm} to the complex analytic orbifold $\fX^+$. Indeed, the real structure $\si:\fX^+\lra\fX^+$ induces a $\C$-antilinear involution $\si$ of the complex vector spaces $\rH^j(\fX^+;\underline{\fL^+})$, in such a way that $\rH^j(\fX;\underline{\fL}) \simeq \Fix_{\si}\, \rH^j(\fX^+;\underline{\fL^+})$, so $\dim_\R \rH^j(\fX;\underline{\fL}) = \dim_\C \rH^j(\fX^+;\underline{\fL^+})$ and $\chi(\fX;\underline{\fL}) = \chi(\fX^+;\underline{\fL^+})$. Setting $\deg\fL:=\deg\fL^+$, one gets: $$\chi(\fX;\underline{\fL}) \ =\ \chi(\fX;O_{\fX}) + \deg(\fL) - 2 \textstyle\sum_{i=1}^k \mathrm{age}_{x_i}(\fL) - \textstyle\sum_{j=1}^{\ell} \mathrm{age}_{y_j}(\fL) $$ where the $(x_i)_{1\leq i\leq k}$ and the $(y_j)_{1\leq j\leq \ell}$ are respectively the cone points and dihedral points of $\fX$. In particular, $\chi(\fX;K_{\fX}^2) = 3(g_{X}-1) +2k+\ell$, where again $g:=\dim\rH^1(\fX;O_{\fX})$ is the genus of $\fX$. % = -3\chi(X) +2k+\ell$. 

\subsection{Spaces of symmetric differentials}\label{computation_of_the_dimension}

As we saw in Section \ref{RR}, if $\fX$ is a compact analytic orbi-curve, then $\chi(\fX;K_{\fX}^2) = \dim\cT(\fX)$ (this is a complex dimension if $\fX$ is defined over $\C$ and a real dimension if $\fX$ is defined over $\R$). While this result is well-known, it is also the $n=2$ case of the \emph{Hitchin parameterisation} of $\Hit(\pifX;\PGL(n;\R))$, as we shall see momentarily. Let us first recall Hitchin's result in the surface group case (\cite{Hitchin_Teich}): If $Y$ is a closed orientable surface of negative Euler characteristic, the choice of a complex analytic structure on $Y$ induces a homeomorphism $$\Hit\big(\piY;\PGL(n;\R)\big)\, \simeq\ \textstyle\bigoplus_{d=2}^n \rH^0(Y;K_Y^d).$$

The main result of \cite{ALS} is the following extension of Hitchin's result to the orbifold case.

\begin{theorem}[\cite{ALS}]\label{Hitchin_param_orb_case}
Let $\fX$ be a compact differentiable orbi-surface of negative Euler characteristic. Then the choice of an analytic structure on $\fX$ induces a homeomorphism $$\Hit\big(\pifX;\PGL(n;\R)\big)\, \simeq\ \textstyle\bigoplus_{d=2}^n \rH^0(\fX;K_{\fX}^d).$$
\end{theorem}

Here, choosing an analytic structure on $\fX$  reduces to choosing a finite Galois cover by a closed orientable surface $Y\lra \fX$ and a complex analytic structure on $Y$ which is preserved by the automorphism group of that cover. As we have seen, the fact that such a cover always exists is a consequence of Selberg's lemma. Note that we are considering at the same time the case where the differentiable orbifold $\fX$ is orientable (so admits a complex analytic structure, i.e.\ the finite group $\pi:=\Aut_X(Y)$ acts holomorphically on $Y$) and the case where it is not (here $Y$ is still a closed orientable surface but $\pi$ will contain orientation-reversing transformation; as a consequence, the coarse moduli space $X$ of $\fX\simeq[Y/\pi]$ will be a differentiable surface with corners that has non-empty boundary or is non-orientable or both).

The proof of Theorem \ref{Hitchin_param_orb_case} consists in adapting Hitchin's proof to the orbifold case. The main tool is the orbifold version of the Non-Abelian Hodge Correspondence (NAHC). In \cite{ALS}, we took a largely equivariant approach to the latter, making the resulting formulation of the NAHC dependent on the choice of a presentation $\fX\simeq [Y/\pi]$. Equivalently, we can rephrase this in terms of $\fG$-Higgs bundles on $\fX$, where $G$ is a real reductive group and $\fG$ is the orbifold group bundle $[(\widetilde{\fX}\times G) / \pifX]$. But in any case, the point is that, if $G$ is the split real form of the adjoint group $\Int(\frg_\C)$, where $\frg_\C$ is a simple complex Lie algebra, then the Hitchin component $\Hit(\piX;G)$ embeds into the moduli space of $\fG$-Higgs bundles, denoted by $\cM_{\fX}(G)$.
\begin{equation}\label{Hitchin_section_diagram}
\xymatrix{
\Hit\big(\pifX,\PGL(n,\R)\big)  \ar@{^{(}->}[r]^{\qquad \tiny{\mathrm{NAHC}}} & \mathcal{M}_{\fX}\big(\PGL(n;\R)\big) \ar[d]^{\mathrm{Hitchin\ fibration}} \\
& \bigoplus_{d=2}^n \rH^0(\fX,K_{\fX}^d) \ar@{.>}@/^/[ul]^{\mathrm{Hitchin\ section\quad}}
}
\end{equation}  In \cite{ALS}, we showed that the \emph{Hitchin fibration}, which is a morphism from the moduli space $\cM_{\fX}(G)$ to a vector space $\mathcal{B}_{\fX}(\frg)$ called the \emph{Hitchin base}, first constructed by Hitchin in the surface group case (\cite{Hitchin_Duke}), was well-defined in the orbifold case. For $\frg=\mathfrak{sl}(n;\R)$, the Hitchin base is $\bigoplus_{d=2}^n \rH^0(\fX;K_{\fX}^d)$, as in Diagram \eqref{Hitchin_section_diagram}. Then we extended Hitchin's construction of a section of that fibration: The image of that section being exactly the embedded copy of $\Hit(\pifX;G)$ in $\cM_{\fX}(G)$, thus proving Theorem \ref{Hitchin_param_orb_case}.

This shows that $\Hit(\pifX;\PGL(n;\R))$ is homeomorphic to the vector space $\mathcal{B}_{\fX}(\mathfrak{sl}(n;\R))$, which is a complex vector space if $\fX$ is complex and a real vector space if $\fX$ is real. Using Theorem \ref{RR_thm}, we can compute the dimension of that vector space. Since we already know how to deduce the result in the real case from the result in the complex case, we will present the proof in the latter case only. From \eqref{degree_powers_of_can_bdle}, we get that, for all $d\in\{2;\,\ldots\,;n\}$, $$\chi(\fX;K_{\fX}^d)\ =\ (2d-1)(g-1) + \textstyle\sum_{i=1}^k R(d,m_i).$$ But for $d\geq 2$, one has $\deg\,K_{\fX}^d=d\deg K_{\fX} > \deg K_{\fX}$, so $\rH^1(\fX;K_{\fX}^d) = 0$ and $$\chi(\fX;K_{\fX})\ =\ \dim\rH^0(\fX;K_{\fX}^d).$$ Thus, when $\fX$ is complex analytic, $\mathcal{B}_{\fX}(\frg)$ is a complex vector space of dimension $$(g-1)\textstyle\sum_{d=2}^n (2d-1) + \textstyle\sum_{d=2}^n \textstyle\sum_{i=1}^k R(d,m_i) = (g-1)(n^2-1) + \textstyle\sum_{d=2}^n \textstyle\sum_{i=1}^k R(d,m_i).$$ The real dimension is twice as much, which indeed coincides with the formula in Theorem \ref{dim_HC_orb_case} (for $\ell=0$).

Let us denote $\PGL(n;\R)$ simply by $G$. A consequence of Theorem \ref{Hitchin_param_orb_case} is that, given an analytic orbi-curve $\fX$, we can embed the Hitchin component $\Hit(\pifX;G)$ into the Hitchin component $\Hit(\pi_1\mathcal{Y};G)$ associated to any Galois cover $\mathcal{Y}\lra\fX$. More precisely, given a Galois cover $\mathcal{Y}\lra\fX$ with automorphism group $\pi$, consider the short exact sequence $$1\lra \pi_1\mathcal{Y} \lra \pifX \lra \pi\lra 1,$$ the induced morphism $\pi \lra \Out(\pi_1\mathcal{Y})$ and the associated action of $\pi$ on $\Hit(\pi_1\mathcal{Y};G)$. Then, the map taking a representation $\rho:\pifX\lra G$ to its restriction $\rho|_{\pi_1\mathcal{Y}}$ induces a homeomorphism $\Hit(\pifX;G)\simeq \Fix_\pi\, \Hit(\pi_1\mathcal{Y};G)$, since $\fX\simeq [\mathcal{Y}/\pi]$ implies that $\rH^0(\fX;K_{\fX}^d) = \Fix_\pi\, \rH^0(\mathcal{Y};K_{\mathcal{Y}}^d)$. As an example of this, consider the Coxeter triangle group $T_{(2,3,7)}$ of \eqref{Coxeter_gp}. It is the orbifold fundamental group of a hyperbolic triangle with vertices of respective orders $2$, $3$ and $7$, which can be obtained as the quotient of the Klein quartic $\cK$ by its full automorphism group. As $\Hit(T_{(2,3,7)};\PGL(6;\R))$ is of (real) dimension $1$ by Theorem \ref{dim_HC_orb_case}, it defines a one-parameter family of Hitchin representations in $\Hit(\pi_1\cK;\PGL(6;\R))$, the latter being, by Hitchin's result for the closed orientable surface $\cK$ (of genus $3$), of real dimension $140$.
\vskip5pt
\textbf{Acknowledgments} It is a pleasure to thank Alexander Schmitt and all participants of the Session on Complex Geometry at ISAAC 2019, as well as Georgios Kydonakis, for feedback on the topics discussed in this paper. I also thank the referee for careful reading and valuable suggestions.

%\bibliographystyle{alpha}
%\bibliography{refs}

\end{document}